\documentclass{amsart}
\begin{document}
\title[Galois groups]{Arithmetic differential equations on $GL_n$, III\\
Galois groups}
\author{Alexandru Buium and Taylor Dupuy}
\def \Rp{R_p}
\def \Rpi{R_{\pi}}
\def \dpi{\d_{\pi}}
\def \bT{{\bf T}}
\def \cI{{\mathcal I}}
\def \cJ{{\mathcal J}}
\def \ZN{\bZ[1/N,\zeta_N]}
\def \tA{\tilde{A}}
\def \o{\omega}
\def \tB{\tilde{B}}
\def \tC{\tilde{C}}
\def \alph{A}
\def \bet{B}
\def \bsigma{\bar{\sigma}}
\def \y{^{\infty}}
\def \Ra{\Rightarrow}
\def \uBS{\overline{BS}}
\def \lBS{\underline{BS}}
\def \lB{\underline{B}}
\def \<{\langle}
\def \>{\rangle}
\def \hL{\hat{L}}
\def \cU{\mathcal U}
\def \cF{\mathcal F}
\def \S{\Sigma}
\def \st{\stackrel}
\def \sd{Spec_{\d}\ }
\def \pd{Proj_{\d}\ }
\def \s{\sigma_2}
\def \i{\sigma_1}
\def \bs{\bigskip}
\def \cD{\mathcal D}
\def \cC{\mathcal C}
\def \cT{\mathcal T}
\def \cK{\mathcal K}
\def \cX{\mathcal X}
\def \sX{X_{set}}
\def \cY{\mathcal Y}
\def \cS{X}
\def \cR{\mathcal R}
\def \cE{\mathcal E}
\def \tcE{\tilde{\mathcal E}}
\def \cP{\mathcal P}
\def \cA{\mathcal A}
\def \cV{\mathcal V}
\def \cM{\mathcal M}
\def \cL{\mathcal L}
\def \cN{\mathcal N}
\def \tcM{\tilde{\mathcal M}}
\def \caS{\mathcal S}
\def \cG{\mathcal G}
\def \cB{\mathcal B}
\def \tG{\tilde{G}}
\def \cF{\mathcal F}
\def \h{\hat{\ }}
\def \hp{\hat{\ }}
\def \tS{\tilde{S}}
\def \tP{\tilde{P}}
\def \tA{\tilde{A}}
\def \tX{\tilde{X}}
\def \tcS{\tilde{X}}
\def \tT{\tilde{T}}
\def \tE{\tilde{E}}
\def \tV{\tilde{V}}
\def \tC{\tilde{C}}
\def \tI{\tilde{I}}
\def \tU{\tilde{U}}
\def \tG{\tilde{G}}
\def \tu{\tilde{u}}
\def \chu{\check{u}}
\def \tx{\tilde{x}}
\def \tL{\tilde{L}}
\def \tY{\tilde{Y}}
\def \d{\delta}
\def \e{\chi}
\def \bW{\mathbb W}
\def \bV{{\mathbb V}}
\def \bF{{\bf F}}
\def \bE{{\bf E}}
\def \bC{{\bf C}}
\def \bO{{\bf O}}
\def \bR{{\bf R}}
\def \bA{{\bf A}}
\def \bB{{\bf B}}
\def \cO{\mathcal O}
\def \ra{\rightarrow}
\def \bx{{\bf x}}
\def \f{{\bf f}}
\def \bX{{\bf X}}
\def \bH{{\bf H}}
\def \bS{{\bf S}}
\def \bF{{\bf F}}
\def \bN{{\bf N}}
\def \bK{{\bf K}}
\def \bE{{\bf E}}
\def \bB{{\bf B}}
\def \bQ{{\bf Q}}
\def \bd{{\bf d}}
\def \bY{{\bf Y}}
\def \bU{{\bf U}}
\def \bL{{\bf L}}
\def \bQ{{\bf Q}}
\def \bP{{\bf P}}
\def \bR{{\bf R}}
\def \bC{{\bf C}}
\def \bD{{\bf D}}
\def \bM{{\bf M}}
\def \bZ{{\mathbb Z}}
\def \xtoleqr{x^{(\leq r)}}
\def \hU{\hat{U}}
\def \k{\kappa}
\def \ee{\overline{p^{\k}}}

\newtheorem{THM}{{\!}}[section]
\newtheorem{THMX}{{\!}}
\renewcommand{\theTHMX}{}
\newtheorem{theorem}{Theorem}[section]
\newtheorem{corollary}[theorem]{Corollary}
\newtheorem{lemma}[theorem]{Lemma}
\newtheorem{proposition}[theorem]{Proposition}
\theoremstyle{definition}
\newtheorem{definition}[theorem]{Definition}
\theoremstyle{remark}
\newtheorem{remark}[theorem]{Remark}
\newtheorem{example}[theorem]{\bf Example}
\numberwithin{equation}{section}
\address{University of New Mexico \\ Albuquerque, NM 87131}
\email{buium@math.unm.edu, taylor.dupuy@gmail.com} 
\maketitle

\begin{abstract}
Differential equations have  arithmetic analogues \cite{book} in which de\-ri\-vatives are replaced by Fermat quotients; these analogues are called arithmetic differential equations and the present paper is concerned with the  ``linear" ones. The equations themselves were introduced in a previous paper \cite{adel2}. In the present paper we deal with the solutions of these equations as well as with the   Galois groups attached to the solutions.
\end{abstract}

\section{Introduction, main definitions, and main results}

In a series of papers beginning with \cite{char} an arithmetic analogue of differential equations was introduced in which derivations are replaced by Fermat quotient operators.  Cf. \cite{book} for an overview. It is then natural to ask for an arithmetic  analogue of linear differential equations. Classically a linear differential equation has  the form 
\begin{equation}
\label{classical}
\frac{d}{dz} U=A\cdot U\end{equation}
where $A$ is, say, a matrix of meromorphic functions on a domain in the complex plane ${\mathbb C}$ with complex variable $z$,   and $U$ is an invertible matrix of unknown meromorphic functions (on a smaller domain). A basic object attached to \ref{classical}  is its differential Galois group which is an algebraic subgroup of $GL_n({\mathbb C})$. This concept is classical, going back to Picard and Vessiot.
A modern version of the theory was developed by Kolchin \cite{kolchin} in the framework of differential algebra. In the present paper we ask for arithmetic analogues, in the spirit of \cite{char, book}, of all of these concepts.
The beginnings of such  a theory   were sketched in \cite{adel2}, where a concept of arithmetic linear differential equation on an algebraic group was introduced; the present paper deals with the solutions of these equations, and especially with the Galois groups attached to these solutions. Our paper  is, in principle, a sequel to \cite{adel1, adel2} but it is entirely independent of these papers. Indeed very little of the theory in \cite{char, book} will be needed here and everything that will be needed will be reviewed in this Introduction. Our main purpose here will be to attach a Galois group to each given solution of a given  linear arithmetic differential equation and to study some basic properties of this group; morally the Galois groups of such equations should (and in some sense will) appear as  subgroups of  ``$GL_n({\mathbb F}_1^a)$" where ${\mathbb F}_1^a$ is the ``algebraic closure of the field with one element"; cf.  \cite{borgerf1} for this interpretation.

\subsection{Main definitions}
We denote by $R$ the unique complete discrete valuation ring  with maximal ideal generated by an odd prime $p$ and with residue field $k=R/pR={\mathbb F}_p^a$, the algebraic closure of ${\mathbb F}_p$. 
   So $R$ can be identified with the ring $W(k)$ of $p$-typical vectors on $k$.
   We denote by 
   $\phi:R\ra R$  the unique ring homomorphism lifting the $p$-power Frobenius on the residue field $k$
   and we denote by
   $\d:R\ra R$ the map $\d a=\frac{\phi(a)-a^p}{p}$.
    We morally view $\d$ as an arithmetic analogue of a derivation \cite{char, book}.
   We denote by 
$R^{\d}$ the monoid of constants $\{\lambda\in R;\d \lambda=0\}$; so $R^{\d}$ consists of $0$ and all roots of unity in $R$. Recall that the reduction mod $p$ map $R^{\d}\ra k$ is an isomorphism of monoids. Also we denote by $K$ the fraction field of $R$. As usual we denote by  ${\mathfrak g}{\mathfrak l}_n(A)$ the ring of $n\times n$ matrices with coefficients in a ring $A$ and we denote by $GL_n(A)$ the group of invertible elements of that ring.  If $A=R$ we will often write 
$$G:=GL_n:=GL_n(R),\ \ {\mathfrak g}:={\mathfrak g}{\mathfrak l}_n:={\mathfrak g}{\mathfrak l}_n(R).$$ 
More generally for a smooth scheme $X$ over $R$ we will often write $X$ instead of $X(R)$. If $u=(u_{ij})\in {\mathfrak g}{\mathfrak l}_n(A)$ then we set
$\phi(u)=(\phi(u_{ij}))$, $\d u=(\d u_{ij})$, $u^{(p)}=(u^p_{ij})$; hence $\phi(u)=u^{(p)}+p\d u$. In what follows we fix a matrix $\Delta(x)\in {\mathfrak g}{\mathfrak l}_n(A)$ with entries in the ring $A=\cO(GL_n)\h=R[x,\det(x)^{-1}]\h$ where $x$ is an $n\times n$ matrix of indeterminates and $\h$  means $p$-adic completion. (This matrix is usually canonically associated to the problem at hand and is uniquely determined by natural symmetry conditions that come with the problem; see \cite{adel2}. We will not be concerned with explaining these conditions  here but rather we will concentrate on the abstract case when $\Delta$ is arbitrary or on specific Examples, cf. \ref{GLn}, \ref{SLn}, \ref{SO(q)} below). Set $\Phi(x)=x^{(p)}+p\Delta(x)$. Moreover for $\alpha\in {\mathfrak g}{\mathfrak l}_n={\mathfrak g}{\mathfrak l}_n(R)$ set 
$$\Delta^{\alpha}(x)=\alpha \cdot \Phi(x)+\Delta(x)=\alpha x^{(p)}+(1+p\alpha)\Delta(x).$$
By
a {\it $\Delta$-linear equation} we will then understand an equation of the form
\begin{equation}
\label{typically}
\d u =\Delta^{\alpha}(u)
\end{equation}
where  $u\in G=GL_n=GL_n(R)$; $u$ is a referred to as a {\it solution} to the equation \ref{typical} and the set 
$G^{\alpha}$ of all $u\in GL_n$ such that \ref{typically} holds
is referred to as the {\it solution set} of \ref{typically}. If we set $\epsilon=1+p\alpha$ and $\Phi^{\alpha}(x)=\epsilon \cdot \Phi(x)$ then \ref{typically} is equivalent to
 \begin{equation}
 \label{dodoly}
 \phi(u)=\Phi^{\alpha}(u).
 \end{equation}
This concept of linearity is always relative to a given $\Delta$. (If $\Delta$ has been fixed and is clear from the context $\Delta$-linear equations are also  referred to as $\d$-linear equations \cite{adel2}.) Note, by the way, that there is a natural concept of equivalence on ${\mathfrak g}{\mathfrak l}_n(A)$ which lies in the background of our discussion; two matrices $\Delta_1$ and $\Delta_2$ in ${\mathfrak g}{\mathfrak l}_n(A)$ are equivalent if and only if there exists $\alpha\in {\mathfrak g}{\mathfrak l}_n(R)$ such that $\Delta_1=\Delta_2^{\alpha}$. We have that $\d u=\Delta_1(u)$ is $\Delta_2$-linear if and only if $\Delta_1$ and $\Delta_2$ are  equivalent.

A function ${\mathcal I}\in R[x,\det(x)^{-1}]\h$ will be called a {\it prime integral} for the $\Delta$-linear equation \ref{typically} if for any solution $u$ of \ref{typically} we have
$$\d({\mathcal I}(u))=0.$$
(Intuitively ${\mathcal I}$ is ``constant" along the solutions of \ref{typically}.) More generally 
an $m$-tuple of functions ${\mathcal I}\in (R[x,\det(x)^{-1}]\h)^m$ is called a prime integral of our equation if each of the components of ${\mathcal I}$ is a prime integral.

The basic examples we have in mind are those in \cite{adel2} and are going to be reviewed below; they are related to the classical groups and for their basic properties we refer to \cite{adel2}. For the purpose of the present article we will not need to review  these properties.

\begin{example}
\label{GLn} We say that $\Delta$ is of type $GL_n$ if
  $\Delta=0$. So in this case $\Phi(x)=x^{(p)}$ and  \ref{typically} and \ref{dodoly} become
\begin{equation}
\label{typical}
\d u =\alpha \cdot u^{(p)}
\end{equation}
and 
 \begin{equation}
 \label{dodo}
 \phi(u)=\epsilon\cdot u^{(p)}
 \end{equation}
 respectively.
It is  worth noting that \ref{dodo} is {\it not} an instance of a linear difference equation in the sense of \cite{SVdP}. Indeed a linear difference equation for $\phi$ has the form 
\begin{equation}
\label{difference}
\phi(u)=\epsilon \cdot u
\end{equation}
rather than the form \ref{dodo}.\end{example}

\begin{example} \label{SLn}
We say that $\Delta$ is of type $SL_n$ if
 $$\Delta(x)=\frac{\lambda(x)-1}{p}\cdot x^{(p)}$$
  where  $p\not| n$ and
\begin{equation}
\label{deflambda}
\lambda(x):=\left( \frac{\det(x^{(p)})}{\det(x)^p}\right)^{-1/n}.
\end{equation}
Here the $-1/n$ power is computed using the usual series  $(1+pt)^a\in \bZ_p[[t]]$ for $a\in \bZ_p$. In this case $\Phi(x)=\lambda(x) \cdot x^{(p)}$ and the equations \ref{typically} and \ref{dodoly} become
\begin{equation}
\label{typical2}
\d u= \left(\lambda(u) \cdot \alpha+\frac{\lambda(u)-1}{p}\right)  \cdot u^{(p)}
\end{equation}
and
\begin{equation}
\label{dodo2}
\phi(u)=\lambda(u)\cdot \epsilon \cdot u^{(p)}
\end{equation}
respectively. Note that, in this case, $\Phi(u)\in SL_n$ for any $u\in SL_n$.
In this context, following \cite{adel2}, it is useful to introduce the $\d$-Lie algebra  ${\mathfrak s}{\mathfrak l}_{n,\d}$ of $SL_n$ as being the set of all $\alpha\in {\mathfrak g}{\mathfrak l}_n$ such that $1+p\alpha\in SL_n$, in other words
$${\mathfrak s}{\mathfrak l}_{n,\d}=\{\alpha\in {\mathfrak g}{\mathfrak l}_n; tr(\alpha)+...+p^{n-1}\det(\alpha)=0\}.$$
This is not a subgroup of $({\mathfrak g}{\mathfrak l}_n,+)$ where $+$ is the usual addition  of matrices but rather a subgroup of $({\mathfrak g}{\mathfrak l}_n,+_{\d})$ where
$a+_{\d}b:=a+b+pab$; the latter group is the group of $R$-points of a group in the category of $p$-adic formal schemes; cf. \cite{adel2}.)
This is in analogy with the Lie algebra ${\mathfrak s}{\mathfrak l}_n$ of $SL_n$ which is given by
$${\mathfrak s}{\mathfrak l}_n=\{\alpha\in {\mathfrak g}{\mathfrak l}_n;tr(\alpha)=0\}.$$
 Note also that if $\alpha\in {\mathfrak s}{\mathfrak l}_{n,\d}$  then ${\mathcal I}(x):=\det(x)$ is a prime integral for the $\Delta$-linear equation $\d u=\Delta^{\alpha}(u)$. Indeed if $u$ is a solution if this equation and 
  $\epsilon=1+p\alpha$  then
$$\begin{array}{rcl}
\phi(\det(u)) & = & \det(\phi(u))\\
\  & = & \det\left(\lambda(u)\cdot  \epsilon \cdot u^{(p)} \right)\\
\  & = & \lambda(u)^n \cdot \det(\epsilon) \cdot \det(u^{(p)})\\
\  & = & \det(u)^p,
\end{array}$$
hence $\d(\det(u))=0$. 
 \end{example}

\begin{example}
\label{SO(q)} 
Let $q\in GL_n$ be defined as
 \begin{equation}
 \label{scorpion3}
 \left(\begin{array}{cl} 0 & 1_r\\-1_r & 0\end{array}\right),\ \ 
\left( 
\begin{array}{ll} 0 & 1_r\\1_r & 0\end{array}\right),\ \ 
\left( \begin{array}{lll} 1 & 0 & 0\\
0 & 0 & 1_r\\
0 & 1_r & 0\end{array}\right),\ \ 
\end{equation}
 where $n=2r, 2r, 2r+1$ respectively. 
Let $SO(q)\subset SL_n$ be the identity component of the subgroup defined by the equations $x^t q x=q$; for $q$ as above $SO(q)$ is denoted by 
$Sp_{2r}, SO_{2r}, SO_{2r+1}$ respectively. We say that $\Delta$ is of type $SO(q)$ if  
$$\Delta(x)=x^{(p)}\cdot \frac{1}{p}(\Lambda(x)-1),$$
 where
$$\Lambda(x)=(((x^{(p)})^tqx^{(p)})^{-1} (x^tqx)^{(p)})^{1/2}.$$
Here, again,  the $1/2$ power is computed using the usual series  $(1+pT)^a\in {\mathfrak g}{\mathfrak l}_n(\bZ_p[[T]])$ for $a\in \bZ_p$, $T=(t_{ij})$. In this case we have $\Phi(x)=x^{(p)}\cdot \Lambda(x)$.
Recall from \cite{adel2} that $\Phi(x)^tq\Phi(x)=(x^tqx)^{(p)}$.
Note also that, in this case, $\Phi(u)\in SO(q)$ for any $u\in SO(q)$; cf. \cite{adel2}.
In this context, following \cite{adel2}, it is useful to introduce the $\d$-Lie algebra  ${\mathfrak s}{\mathfrak o}(q)_{\d}$ of $SO(q)$ as being the set of all $\alpha\in {\mathfrak g}{\mathfrak l}_n$ such that $1+p\alpha\in SO(q)$, in other words
$${\mathfrak s}{\mathfrak o}(q)_{\d}=\{\alpha\in {\mathfrak g}{\mathfrak l}_n; \alpha^t q+q\alpha+p\alpha^tq\alpha=0\}.$$
This is, again, a subgroup of $({\mathfrak g}{\mathfrak l}_n,+_{\d})$;
and this is, again,  in analogy with the Lie algebra ${\mathfrak s}{\mathfrak o}(q)$ of $SO(q)$ which is given by
$${\mathfrak s}{\mathfrak o}(q)=\{\alpha\in {\mathfrak g}{\mathfrak l}_n;\alpha^t q+q\alpha=0\}.$$
Note also that if $\alpha\in {\mathfrak s}{\mathfrak o}(q)_{\d}$ then ${\mathcal I}(x):=x^tqx$ is a prime integral for
the $\Delta$-linear equation $\d u=\Delta^{\alpha}(u)$. Indeed, if $u$ is a solution of this equation and $\epsilon=1+p\alpha$ then, using the identity $\Phi(x)^t q \Phi(x)=(x^t qx)^{(p)}$,  we get
$$\begin{array}{rcl}
\phi(u^tqu) & = & \phi(u)^t q \phi(u)\\
\  & = & \Phi(u)^t \epsilon^t q \epsilon \Phi(u)\\
\  & = & \Phi(u)^t q  \Phi(u)\\
\  & = & (u^tqu)^{(p)},
\end{array}
$$
which  implies $\d(u^tqu)=0$.
\end{example}

\subsection{Main results}
 One has an existence and uniqueness result for our equations \ref{typically}; cf.  Propositions \ref{exu}, \ref{mama}, \ref{tata}, and Remark \ref{eggs}  in the body of the paper: 

\begin{theorem}\label{drink}
Let $u_0\in GL_n$ and $\alpha\in {\mathfrak g}{\mathfrak l}_n$ and let $\Delta$ be arbitrary. Then the following hold:

1) There is a unique  $u\in GL_n$ satisfying \ref{typically} such that $u\equiv u_0$ mod $p$.

 2) If $\Delta$, $u_0$, and $\alpha$  have entries in a complete valuation subring $\cO$ of $R$ then  $u$ also has entries in $\cO$.
 
 3) If $u_0\in SL_n$,  $\alpha\in {\mathfrak s}{\mathfrak l}_{n,\d}$, and $\Delta$ is of type $SL_n$ then $u\in SL_n$.
 
 4) If $u_0\in SO(q)$,  $\alpha\in {\mathfrak s}{\mathfrak o}(q)_{\d}$, and $\Delta$ is of type $SO(q)$ then $u\in SO(q)$.

5) If $u_0$ and $\alpha$  have entries in a valuation $\d$-subring $\cO$ of $R$ with  
finite residue field and either $\Delta$ is of type $GL_n$ (i.e. $\Delta=0$) or
$\Delta$ is of type $SL_n$ and $u\in SL_n$
 then $u$ has entries in a $\d$-subring of $R$ which is generically finite over $\cO$.
\end{theorem}

Here by a $\d$-subring $\cO$ of $R$ we understand a subring with $\d \cO\subset \cO$. By a valuation subring of $R$ we mean the intersection of $R$ with a subfield of the field of fractions $K$ of $R$. Also an extension of integral domains is called generically finite if the  induced extension between fraction fields is  finite. 

The above theorem allows us to introduce the first steps in a $\d$-Galois theory attached to $\Delta$-linear equations \ref{typical}. In particular we will attach $\d$-Galois groups to such equations and prove results about their form in ``generic" cases. Here are some details. Start with a  $\d$-subring $\cO\subset R$, let $\alpha\in {\mathfrak g}{\mathfrak l}_n(\cO)$ and let $u\in GL_n(R)$
be a solution of \ref{typically}. Consider the subring $\cO\{u\}$ of $R$ generated by $\cO$ and $u,\d u, \d^2 u,...$; so $\d \cO\subset \cO$. Consider  the group  $Aut_{\d}(\cO\{u\}/\cO)$  of all $\cO$-algebra automorphisms $\sigma$ of $\cO\{u\}$ such that $\sigma\circ \d=\d\circ \sigma$ on $\cO\{u\}$. Consider furthermore the subgroup
$\tilde{G}_{u/\cO}$ of $Aut_{\d}(\cO\{u\}/\cO)$ consisting of all $\sigma\in Aut_{\d}(\cO\{u\}/\cO)$ such that $u^{-1}\cdot \sigma(u)\in GL_n(\cO)$. There is natural map, which is an injective group homomorphism,
 \begin{equation}
 \label{torpedo}
 \tilde{G}_{u/\cO}\ra GL_n(\cO)\end{equation}
sending any $\sigma$ into $c_{\sigma}:=u^{-1}\cdot \sigma(u)$. Finally define the $\d$-{\it Galois group} of $u/\cO$ as the image $G_{u/\cO}$ of  \ref{torpedo}. In particular $G_{u/\cO}\simeq \tilde{G}_{u/\cO}$.
Our next task is to ``compute/bound" $\d$-Galois groups.
We begin with $\Delta$ of type $SL_n$ and $SO(q)$:

\begin{theorem}\label{shesleeps}
\

1) Assume $\Delta$ is of type $SL_n$ and  let $\alpha\in {\mathfrak g}{\mathfrak l}_n(\cO)$,  $u\in G^{\alpha}$. Then for any $c\in G_{u/\cO}$
we have $\d(\det(c))=0$.

2) Assume $\Delta$ is of type $SO(q)$ and  let $\alpha\in {\mathfrak g}{\mathfrak l}_n(\cO)$, $u\in SO(q)\cap G^{\alpha}$. Then for any $c\in G_{u/\cO}$ we have $\d(c^t q c)=0$.
\end{theorem}

Cf.  Propositions \ref{forget} plus \ref{clock}.

Theorem \ref{shesleeps} shows that if $\Delta$ is of type $SL_n$ or $SO(q)$ the $\d$-Galois group $G_{u/\cO}$ is ``close to being contained" in $SL_n$ and $SO(q)$ respectively (provided $u$ is in these groups respectively). Indeed 
$G_{u/\cO}$ being contained in $SL_n$ (respectively $SO(q)$)  means 
 $\det(c)=1$ (respectively $\det(c)=1$ and $c^tqc=q$) for $c\in G_{u/\cO}$.  Theorem \ref{shesleeps}, however, 
merely guarantees that $\d(\det(c))=0$ or $\d(c^tqc)=0$, which is a ``slightly" weaker property.

Next we consider the case $\Delta$ is of type $GL_n$ (i.e. $\Delta=0$).
To state our result below we let $W\subset G$ be the Weyl group of all matrices obtained from the identity matrix by permuting its columns. Let $T\subset G$ be the maximal torus of diagonal matrices with entries in $R$ and consider
the normalizer $N=WT=TW$ of $T$ in $G$. We denote by $1\in G$ the identity matrix.
Also consider the subset (not a subgroup!) $G^{\d}$ of $G$ consisting of all elements of $G$ with entries in the monoid of constants  $R^{\d}$. Let $N^{\d}=N\cap G^{\d}$ and $T^{\d}=T\cap G^{\d}$.
Then $N^{\d}$ and $T^{\d}$ are  subgroups (not just subsets!) of $G$. Also
$N^{\d}=WT^{\d}=T^{\d}W$.
We also use below the notation $K^a$ for the algebraic closure of the fraction field $K$ of $R$; the Zariski closed sets  of $GL_n(K^a)$ are then viewed as (possibly reducible) varieties over $K^a$. A subgroup of $GL_n(K^a)$ is called diagonalizable if it is conjugate in $GL_n(K^a)$ to a subgroup of the group of diagonal matrices.
The next result illustrates some ``generic" features of our $\d$-Galois groups in case $\Delta=0$; assertion 1) 
shows that the $\d$-Galois
group is generically ``not too large". Assertion 2) 
 shows that the $\d$-Galois groups are generically ``as large as possible". As we shall see presently, the meaning of the word {\it generic} is different in each of the $2$ situations: in situation 1) {\it generic} means {\it outside a Zariski closed set};  in situation 2)  {\it generic} means {\it outside a set of the first category} (in the sense of Baire category).

\begin{theorem}
\label{food} Assume $\Delta(x)=0$.

1) There exists a Zariski closed subset $\Omega\subset GL_n(K^a)$ not containing $1$ such that for any $u\in GL_n(R) \backslash \Omega$ the following holds.  Let $\alpha=\d u \cdot (u^{(p)})^{-1}$ and let $\cO$ be a valuation $\d$-subring of $R$ containing $\alpha$. Then $G_{u/\cO}$ contains a normal subgroup of finite index which is diagonalizable.  

2) There exists  a subset $\Omega$  of the first category in the metric space
$$X=\{u\in GL_n(R);u\equiv 1\ \ \text{mod}\ \ p\}$$
such that for any $u \in X\backslash \Omega$ the following holds.
 Let $\alpha=\d u\cdot (u^{(p)})^{-1}$. Then there exists a  valuation $\d$-subring $\cO$ of $R$ containing $R^{\d}$ such that
$\alpha\in {\mathfrak g}{\mathfrak l}_n(\cO)$ and such that $G_{u/\cO}=N^{\d}$.

\end{theorem}

Cf.   Propositions \ref{dimless},  \ref{transcendence},   in the body of the paper.

The groups $W$ and  $N^{\d}$ should be morally viewed as ``incarnations" of   ``$GL_n({\mathbb F}_1)$" and ``$GL_n({\mathbb F}_1^a)$" where 
 ``${\mathbb F}_1$" and ``${\mathbb F}_1^a$" are the ``field with element" and ``its algebraic closure" respectively; cf.  \cite{borgerf1}.  This suggests that the $\d$-Galois theory we are proposing here should be viewed as a Galois theory over ``${\mathbb F}_1$". By the way  Theorem \ref{food} suggests the following question: {\it Is the $\d$-Galois group  $G_{u/\cO}$  {\it always} a subgroup of $N$?}
 The answer to this turns out to be negative in general (cf. Example \ref{counterexample}) but something close to an affirmative answer may still be true.
 
 \medskip

We end with a couple of remarks comparing the theory above  with some familiar situations.

\begin{remark}\label{expostuff}
It is worth comparing Equation
 \ref{dodo} with the familiar linear equations in analysis in the case  $n=1$; in case $n=1$ Equation \ref{dodo} is, of course, 
\begin{equation}\label{dododo}
\phi(u)=\epsilon \cdot u^p
\end{equation}
where $\epsilon=1+p\alpha$, $\alpha\in R$, $u\in R^{\times}$. This equation can be solved as follows. Write $\epsilon=\exp(p\beta)$, where $\exp:pR\ra 1+pR$ is the group isomorphism given by the $p$-dic exponential and $\beta\in R$. Then the set of solutions to \ref{dododo} consists of all $u\in R^{\times}$ of the form
\begin{equation}
\label{research}
u=\zeta \cdot \exp\left( \sum_{n=1}^{\infty} p^n \phi^{-n}(\beta)\right)
\end{equation}
where $\zeta\in R^{\times}$, $\d \zeta=0$. On the other hand consider the group homomorphism $\psi:R^{\times}\ra R$ defined by
\begin{equation}
\label{psimap}
u\mapsto \psi(u)=\frac{1}{p}\log\left(\frac{\phi(u)}{u^p}\right)=\sum_{n=1}^{\infty}(-1)^{n-1}
\frac{p^{n-1}}{n}\left(\frac{\d u}{u^p}\right)^n
\end{equation}
where $\log$ is the $p$-adic logarithm.
Then Equation \ref{dododo} is equivalent to the equation
\begin{equation}
\label{teach}
\psi(u)=\beta
\end{equation}
Now the homomorphism $\psi$ above should be viewed as an analogue
of the logarithmic derivative map ${\mathcal M}(D)^{\times}\ra {\mathcal M}(D)$, $$u\mapsto u'/u,$$ where ${\mathcal M}(D)$ is the field of meromorphic functions on a disk $D\subset {\mathbb C}$, say, and $u'=\frac{du}{dz}$, where $z$ is a complex variable. So the analogue, in analysis,  of Equation \ref{teach} is the equation 
\begin{equation}
\label{dududu}
\frac{u'}{u}=\beta,
\end{equation}
where $\beta \in {\mathcal M}(D)$. For $\beta$ holomorphic in $D$ the solutions to Equation \ref{dududu} are of the form
\begin{equation}
\label{junior}
u=c\cdot \exp\left(\int \beta dz\right)
\end{equation}
where $\exp$ is the complex exponential and $c\in {\mathbb C}$. Hence the elements \ref{research} in $R^{\times}$
should be viewed as  arithmetic analogues of the functions \ref{junior} in ${\mathcal M}(D)$. \end{remark}

\begin{remark}
It is  worth comparing the $\Delta$-linear equations \ref{typical} with Lang's framework in \cite{lang}. Indeed in \cite{lang} Lang considers the map
\begin{equation}
\label{characteristic}
GL_n(k)\ra GL_n(k),\ \ a\mapsto a^{(p)}\cdot a^{-1},
\end{equation}
where $k$ is an algebraically closed field of characteristic $p$. 
This is a non-abelian cocycle for the adjoint action of $GL_n(k)$ on itself.
A natural lift of \ref{characteristic} to characteristic zero is the map
\begin{equation}
\label{tv}
GL_n(R)\ra GL_n(R),\ \ a\mapsto \phi(a)\cdot a^{-1}.\end{equation}
The fiber of \ref{tv}  over $\alpha\in {\mathfrak g}{\mathfrak l}_n(R)$ consists of the solutions $u\in GL_n(R)$ to the  {\it linear difference equation} \ref{difference} which, as already noted, 
is quite different from the equation \ref{dodo}. By the way the  equation \ref{difference} can be studied in at least two ways leading to two rather different theories: one way is from the viewpoint of  difference algebra \cite{SVdP}; the other way is from the  $\d$-arithmetic viewpoint \cite{book}. The 
$\d$-arithmetic viewpoint on  equations \ref{difference} tends to lead to profinite groups; our $\d$-arithmetic study  of the equations \ref{dodo} will lead to torsion  groups (hence to inductive, rather than projective, limits of finite groups). This makes the $\d$-arithmetic study of equations \ref{dodo} and the $\d$-arithmetic study of equations \ref{difference} quite different in nature. 
 Neverthless there are cases (such as that of abelian varieties \cite{char}) where one encounters combinations of profinite and torsion groups; so it is conceivable that the $\d$-arithmetic theories of \ref{dodo} and \ref{difference} can be unified.

On the other hand 
\ref{characteristic} has another natural lift to characteristic zero which is
\begin{equation}
\label{coveringg}
GL_n(R)\ra GL_n(R),\ \ a\mapsto a^{(p)}\cdot a^{-1}.\end{equation}
(This map is not induced by an endomorphism of the scheme  $GL_n$ but rather by an endomorphism of the $p$-adic completion of $GL_n$.)
Composing this with inversion $b\mapsto b^{-1}$ one gets a map
\begin{equation}
\label{covering}
GL_n(R)\ra GL_n(R),\ \ a\mapsto a\cdot (a^{(p)})^{-1}.\end{equation}
Note now that the set of solutions to any of the  equations \ref{dodo} is a fiber of the map
\begin{equation}
\label{coverin}
GL_n(R)\ra GL_n(R),\ \ a\mapsto \phi(a)\cdot (a^{(p)})^{-1}\end{equation}
But \ref{covering} and \ref{coverin} induce by restriction the same map $GL_n(\bZ_p)\ra GL_n(\bZ_p)$. This connection points towards a link between the arithmetic of usual coverings such as  \ref{coveringg} and the ``$\d$-Galois theory" of $\Delta$-linear equations such as \ref{dodo}. Also, in some sense, our paradigm here can be viewed as a lift to characteristic zero, in the framework of ``$\d$-geometry", of Lang's  characteristic $p$ algebro-geometric paradigm. \end{remark}

The paper is organized as follows. In section 2 we provide the proof of Theorem \ref{drink}. In section 3 we amplify our definitions and foundational discussion and we prove, in particular,  Theorems \ref{shesleeps} and  \ref{food}.

\medskip

\subsection{Acknowledgement}
The authors are indebted to P. Cartier for inspiring discussions. Also the first author would like to acknowledge partial support from the Hausdorff Institute of Mathematics in Bonn, from the NSF through grant DMS 0852591, 
 from the Simons Foundation
(award 311773), 
and from the Romanian National Authority
for Scientific Research, CNCS - UEFISCDI, project number
PN-II-ID-PCE-2012-4-0201.

\section{Existence, uniqueness, and rationality of solutions}

The following proposition  is an existence and uniqueness result for solutions of $\Delta$-linear equations. In the Propositions below $\Delta(x)$ is arbitrary unless otherwise stated and, as usual, $\Phi(x)=x^{(p)}+p\Delta(x)$.

\begin{proposition}
\label{exu}
Let  $u_0\in GL_n(R)$, and $\alpha \in {\mathfrak g}{\mathfrak l}_n(R)$. Then the $\Delta$-linear equation 
$\d u=\alpha \cdot \Phi(u)+ \Delta(u)$ has a unique solution $u\in GL_n(R)$ such that $u\equiv u_0$ mod $p$. 
\end{proposition}

{\it Proof}.
Recall that the equation above is equivalent to $\phi(u)=\epsilon  \cdot \Phi(u)$ where $\epsilon=1+p\alpha$. 
To check the uniqueness of the solution assume $\phi(u)=\epsilon \cdot \Phi(u)$ and
$\phi(v)=\epsilon \cdot \Phi(v)$ with $u,v \in GL_n(R)$, $u\equiv v$ mod $p$. Then we prove by induction that $u\equiv v$ mod $p^n$. Indeed if the latter is the case then $u^{(p)}\equiv v^{(p)}$ mod $p^{n+1}$ and $\Delta(u)\equiv \Delta(v)$ mod $p^{n+1}$ hence $\Phi(u)\equiv \Phi(v)$ mod $p^{n+1}$. Hence $\phi(u)\equiv \phi(v)$ mod $p^{n+1}$. Hence $u\equiv v$ mod $p^{n+1}$. 

To check the existence of a solution $u$ such that $u\equiv u_0$ mod $p$ we define a sequence of matrices $u_n\in GL_n(R)$ by the formula
$$u_{n+1}=\phi^{-1}(\epsilon \cdot \Phi(u_n)),\ \ n\geq 0.$$
We claim that for all $n\geq 0$ we have 
$$\phi(u_n)\equiv  \epsilon \cdot \Phi(u_n)\ \ \ \text{mod}\ \ p^{n+1}.$$
Assuming the claim we get $u_{n+1}\equiv u_n$ mod $p^{n+1}$ hence $u_n$ converges $p$-adically to some $u\in GL_n(R)$.
Also $\phi(u)=\epsilon \cdot \Phi(u)$ which ends our proof. We are left with checking the claim. We proceed by induction. The case $n=0$ is clear. Assume now $\phi(u_n)  \equiv \epsilon \cdot \Phi(u_n)$ mod $p^{n+1}$. Hence
$$\phi^{-1}(\epsilon \cdot \Phi(u_n))\equiv u_n\ \ mod\ \ p^{n+1},$$
hence
$$\Phi(\phi^{-1}(\epsilon \cdot \Phi(u_n)))\equiv \Phi(u_n)\ \ mod\ \ p^{n+2}.$$
Hence
$$\begin{array}{rcl}
\epsilon\cdot \Phi(u_{n+1}) & = & \epsilon\cdot \Phi(\phi^{-1}(\phi(u_{n+1})))\\
\  & \ & \ \\
\  & = & \epsilon \cdot \Phi(\phi^{-1}(\epsilon\cdot \Phi(u_n)))\\
\  & \ & \ \\
\  & \equiv & \epsilon \cdot \Phi(u_n)\ \ \ \text{mod}\ \ \ p^{n+2}\\
\  & \ & \ \\
\  & = & \phi(u_{n+1}),\end{array}
$$
and the induction step follows. \qed

\begin{remark}
If, in Proposition \ref{exu},  $\Delta=0$, $n=1$, and $u_0\equiv \zeta$ mod $p$ where $\zeta\in R$ is a root of unity, the solution $u$  has a closed  form:
$$u=\zeta \cdot \epsilon_{-1} \cdot \epsilon_{-2}^p\cdot \epsilon_{-3}^{p^2}\cdot ...\ \ \text{(convergent product)}$$
where $\epsilon_i=\phi^i(\epsilon)$ for $i\in \bZ$. This computation implies the formula in Remark \ref{expostuff}.
\end{remark}

\begin{remark}
\label{eggs}
If in Proposition \ref{exu} we have $\Delta$ of type $SL_n$, $u_0\in SL_n(R)$, and $\alpha\in {\mathfrak s}{\mathfrak l}_{n,\d}$ then $u\in SL_n(R)$. Indeed this follows because $\Phi(a)\in SL_n(R)$ and $\phi^{-1}(a)\in SL_n(R)$ for all $a\in SL_n(R)$; hence if $u_n$ is as in the proof of that Proposition 
then $u_n\in SL_n(R)$. Similarly if $\Delta$ is of type $SO(q)$, $u_0\in SO(q)$, and $\alpha\in {\mathfrak s}{\mathfrak o}(q)_{\d}$ then $u\in SO(q)$. The above proves assertions 3 and 4 in Theorem \ref{drink}.
Alternatively these assertions can be deduced as follows. Let ${\mathcal I}(x)$ be $\det(x)$ or $x^tqx$ respectively and let $u$ be such that
$\d u=\Delta^{\alpha}(u)$ with $u\equiv u_0$ mod $p$ with either $u_0\in SL_n$, $\alpha\in {\mathfrak s}{\mathfrak l}_{n,\d}$ or
$u_0\in SO(q)$, $\alpha\in {\mathfrak s}{\mathfrak o}(q)_{\d}$ respectively. By the discussion in Examples \ref{SLn} and \ref{SO(q)} we have $\d({\mathcal I}(u))=0$ hence ${\mathcal I}(u)$ is either $0$ or a root of unity in $R$.
 On the other hand we have ${\mathcal I}(u)\equiv {\mathcal I}(u_0)$ mod $p$ hence, since ${\mathcal I}(u_0)=0$, 
  we have ${\mathcal I}(u)\equiv 0$ mod $p$. 
 Since ${\mathcal I}(u)$ is either $0$ or a root of unity we conclude it must be $0$, hence $u\in SL_n$ or $u\in SO(q)$ respectively.
\end{remark}

\begin{remark}
In notation of Propositon \ref{exu} the natural reduction map $G^{\alpha}\ra GL_n(k)$ is a bijection. So each solution set $G^{\alpha}$ has a natural structure of group; but of course with this structure $G^{\alpha}$ is not a subgroup of $GL_n(R)$.
\end{remark}

Let us address the question of ``rationality" of solutions of $\Delta$-linear equations.

Let $\cO\subset R$ be a subring. 
Recall that $\cO$ is called a 
$\d$-subring if 
$\d\cO\subset \cO$. Also we say $\cO$ is a
a valuation subring of $R$ if $\cO$ is the intersection of $R$ with a subfield of $K$. Any valuation subring of $R$ is a discrete valuation ring with maximal ideal generated by $p$.
 Note that if $\cO$ is a valuation subring which is complete then either  $\cO=R$ or there exists $\nu\geq 1$ such that $\cO=R^{\phi^{\nu}}$, the fixed ring of $\phi^{\nu}$; in particular such an $\cO$ is automatically a $\d$-subring. An extension $\cO\subset \cO'$ of  subrings of $R$ will be called generically finite if the extension of their fraction fields is finite; if in addition $\cO$ is a valuation subring then  $\cO'$ is  a localization  of a finite extension of $\cO$; if, in addition $\cO$ is complete then any generically finite extension of $\cO$ in $R$ is finite.

\begin{proposition}
\label{mama}
Assume $\cO$ is a complete valuation subring of $R$ (hence also a $\d$-subring).
If in Proposition \ref{exu} we have 
$$\Delta\in {\mathfrak g}{\mathfrak l}_n(\cO[x,\det(x)]\h),\ \ u_0\in GL_n(\cO),\ \ \alpha\in {\mathfrak g}{\mathfrak l}_n(\cO)$$ then 
$u\in GL_n(\cO)$.
\end{proposition}

{\it Proof}.
Let $\cO=R^{\phi^{\nu}}$. Then $\phi^{\nu}(u_0)=u_0$ and $\phi^{\nu}(\alpha)=\alpha$ hence $\phi^{\nu}(\epsilon)=\epsilon$, where $\epsilon=1+p\alpha$. Also $\phi^{\nu}(\Delta(a))=\Delta(\phi^{\nu}(a))$, and hence $\phi^{\nu}(\Phi(a))=\Phi(\phi^{\nu}(a))$, for all $a\in GL_n(R)$.
Since $\phi(u)=\epsilon \cdot \Phi(u)$ and $u\equiv u_0$ mod $p$ it follows that
$$\phi^{\nu+1}(u)=\phi^{\nu}(\epsilon) (\phi^{\nu}(\Phi(u)))=\epsilon \cdot \Phi( (\phi^{\nu}(u)))$$
 and $\phi^{\nu}(u)\equiv \phi^{\nu}(u_0)\equiv u_0$ mod $p$. By the uniqueness in Proposition \ref{exu} it follows that $\phi^{\nu}(u)=u$ hence $u\in GL_n(\cO)$.
\qed

\begin{proposition}
\label{tata}
Assume $\cO$ is a valuation $\d$-subring of $R$ with finite residue field. 
Assume in Proposition \ref{exu} that one of the following holds:

1) $\Delta$ is of type $GL_n$ (i.e. $\Delta=0$), $u_0\in GL_n(\cO)$, and $\alpha\in {\mathfrak g}{\mathfrak l}_n(\cO)$.

2) $\Delta$ is of type $SL_n$, $u_0\in SL_n(\cO)$, and $\alpha\in {\mathfrak s}{\mathfrak l}_{n,\d}\cap {\mathfrak g}{\mathfrak l}_n(\cO)$.

Then there exists a generically  finite extension of $\d$-subrings $\cO\subset \cO'$ of $R$ such that
$u\in GL_n(\cO')$.
\end{proposition}

{\it Proof}. Assume we are in case 2; case 1 is similar (and indeed slightly easier).

By Proposition \ref{mama} if $\widehat{\cO}$ is the completion of $\cO$ then $u\in GL_n(\widehat{\cO})$ hence there exists $\nu\geq 0$ such that
 $\phi^{\nu+1}(u)=u$. Let  $N=n^2$ and identify 
 the points of ${\mathbb A}^N$ with $n\times n$ matrices.
 Let
 $$\lambda_{\nu}(u)=\phi^{\nu}(\lambda(u))\cdot \phi^{\nu-1}(\lambda(u))^p\cdot ...\cdot
 \lambda(u)^{p^{\nu}}.$$
 Using $\phi(u)=\lambda(u) \cdot \epsilon \cdot u^{(p)}$, and setting $\epsilon_j=\phi^j(\epsilon)$, we get 
\begin{equation}
\label{marga}
u=\phi^{\nu+1}(u)=\lambda_{\nu}(u)\cdot \varphi(u),\end{equation}
where $\varphi:{\mathbb A}^{N}\ra {\mathbb A}^{N}$ is the morphism  of schemes over  $\cO$ defined on points by
$$\varphi(v)=\epsilon_{\nu}(\epsilon_{\nu-1}(\epsilon_{\nu-2}(...(\epsilon v^{(p)})^{(p)})^{(p)}...)^{(p)}.$$
Let $K^a$ be an algebraic closure of $K$, let $F$ be the fraction field of $\cO$,  and let $F^a$ be the algebraic closure of $F$ in $K^a$. Note that $\varphi:{\mathbb A}^N(K^a)\ra {\mathbb A}^N(K^a)$  is obtained by composing maps $\eta\mapsto \epsilon_j \eta$ with copies of the  map $\eta\mapsto \eta^{(p)}$;
both these types of maps are given by homogeneous polynomials (of degree $1$ and $p$ respectively) and have the property that the pre-image of $0$ is $0$. Hence $\varphi$ is given by
$$\varphi(\eta)=(\Phi_1(\eta),...,\Phi_N(\eta))$$
where $\Phi_1,...,\Phi_N\in F[x_1,...,x_N]$ are homogeneous polynomials of  degree $p^{\nu+1}>1$ and $\varphi^{-1}(0)=\{0\}$; hence $\Phi_1,...,\Phi_N$ have no common zero in ${\mathbb A}^N(K^a)$ except at the origin.  Consider an extra variable $x_0$ and  consider the projective variety $V\subset {\mathbb P}^{N}$ defined by the equations 
\begin{equation}
\label{cheese}
\Phi_j(x_1,...,x_N)-x_0^{p^{\nu+1}-1}x_j=0.\end{equation}
Clearly the intersection of $V$ with the hyperplane $x_0=0$ is empty. So $V$ has dimension zero hence  $V(K^a)$ is finite. Since $V$ is defined over $F$ we have
$V(K^a)=V(F^a)$.
By equation \ref{marga} the point
$$(\lambda_{\nu}(u)^{-1/(p^{\nu+1}-1)}:u)\in {\mathbb P}^{N}(K)$$
belongs to $V(K)$ hence it belongs to $V(F^a)$. (Here the $1/(p^{\nu+1}-1)$-power is computed, again, using the series $(1+pt)^a\in \bZ_p[[t]]$ for $a\in \bZ_p^{\times}$). It follows that 
\begin{equation}
\label{ham}
\lambda_{\nu}(u)^{1/(p^{\nu+1}-1)}\cdot u\in {\mathbb A}^{N}(F^a)={\mathfrak g}{\mathfrak l}_n(F^a)\end{equation}
hence 
$$\det(\lambda_{\nu}(u)^{1/(p^{\nu+1}-1)}\cdot u)\in F^a.$$
Since, by Remark \ref{eggs}, $\det(u)=1$ we get
$(\lambda_{\nu}(u)^{1/(p^{\nu+1}-1)})^n\in F^a$ hence 
$$\lambda_{\nu}(u)^{1/(p^{\nu+1}-1)}\in F^a.$$
By \ref{ham} again we get $u\in GL_n(F^a)$ which ends the proof.
\qed

\bigskip

Note that the Propositions in this section imply Theorem \ref{drink} in the Introduction. The consideration of the variety cut out by equations \ref{cheese} is a trick from \cite{FS} and is  an indication of an interesting  link between the paradigm of the present paper and the arithmetic of dynamical systems on projective space.


\section{$\d$-Galois groups}

Recall that $\d$-Galois groups were defined in the Introduction. We will review  here the notation involved and some related concepts. Then we will prove  a series of Propositions amounting to Theorem \ref{food}.

As usual we often denote by $G$ the group $GL_n(R)$ and by ${\mathfrak g}{\mathfrak l}_n$ the Lie algebra ${\mathfrak g}{\mathfrak l}_n(R)$.
Let $\Delta(x)\in {\mathfrak g}{\mathfrak l}_n(R[x,\det(x)^{-1}]\h)$, $x$  an $n\times n$ matrix of indeterminates,
and let $\Phi(x)=x^{(p)}+p\Delta(x)$. 
Let  $\alpha \in {\mathfrak g}{\mathfrak l}_n$, $\Delta^{\alpha}(x)=\alpha\cdot \Phi(x)+\Delta(x)$,
 and consider the $\Delta$-linear equation 
\begin{equation}
\label{tree1}
\d u=\Delta^{\alpha}(u).\end{equation}
Recall that if $\Phi^{\alpha}(x)=\epsilon\cdot \Phi(x)$, $\epsilon=1+p\alpha$, then this equation is equivalent to the equation 
\begin{equation}
\label{tree2}
\phi(u)=\Phi^{\alpha}(u).\end{equation}
 Let $G^{\alpha}$ be the set of solutions to Equation \ref{tree1}, let $u\in G^{\alpha}$ be a fixed solution, let $\Phi_u(x)=\Phi(u)^{-1}\Phi(ux)$, $\Delta_u(x)=\frac{1}{p}(\Phi_u(v)-v^{(p)})$,
 and let $G_u$ be the set of solutions $v\in G$ to the $\Delta_u$-linear equation 
\begin{equation}
\label{tree3}
\d v=\Delta_u(v),\end{equation}
equivalently to the equation 
\begin{equation}
\label{cry}
\phi(v)=\Phi_u(v).\end{equation}
Note that
$$uG_u\subset G^{\alpha}.$$ Indeed if $c\in G_u$ we have
$$\phi(uc)=\phi(u)\cdot \phi(c)=\epsilon \cdot \Phi(u)\cdot \phi(c)=\epsilon\cdot \Phi(uc)$$
so $uc\in G^{\alpha}$.

Let now $\cO$ be a $\d$-subring of $R$. Assume $\alpha\in {\mathfrak g}{\mathfrak l}_n(\cO)$ and let $u\in GL_n(R)$
be a solution of Equation \ref{tree1}. Recall from the Introduction 
 the group  $Aut_{\d}(\cO\{u\}/\cO)$  of all $\cO$-algebra automorphisms $\sigma$ of $\cO\{u\}$ such that $\sigma\circ \d=\d\circ \sigma$ on $\cO\{u\}$, its  subgroup $\tilde{G}_{u/\cO}$ and the injective group homomorphism  $\tilde{G}_{u/\cO}\ra GL_n(\cO)$
sending any $\sigma$ into $c_{\sigma}:=u^{-1}\cdot \sigma(u)$. Then the $\d$-Galois group $G_{u/\cO}$ was defined as the image $G_{u/\cO}$ of the homomorphism
 $\tilde{G}_{u/\cO}\ra GL_n(\cO)$.
 
 In the special cases of interest to us the $\d$-Galois group has
 a ``$\d$-theoretic description/bound" which we now discuss.
Let $x',x'',...$ be new matrices of indeterminates and consider the polynomial ring
$\cO\{x\}:=\cO[x,x',x'',...]$.
There is a unique ring endomorphism $\phi$ of $\cO\{x\}$ whose restriction to $\cO$ is $\phi$ and such that $\phi(x)=x^{(p)}+px'$, $\phi(x')=(x')^{(p)}+px''$, etc. Define the map $\d:\cO\{x\}\ra \cO\{x\}$ by $\d f=p^{-1}(\phi(f)-f^p)$.
We let $I_{u/\cO}$ be the kernel of the unique $\cO$-algebra map $\cO\{x\}\ra R$, sending $x\mapsto u$, $x'\ra \d u$, $x''\mapsto \d^2 u$, etc. (the ideal of $\d$-algebraic relations among the entries of $u$); note that 
$\cO\{u\}$ is then the image of the map $\cO\{x\}\ra R$ above. We let
$\Sigma_{u/\cO}$ be the subgroup of $GL_n(\cO)$ consisting of all matrices $c$ such that the $\cO$-automorphism $\sigma_c:\cO\{x\}\ra \cO\{x\}$ defined by $\sigma_c(x)=xc$, $\sigma(x')=\d(xc)$, $\sigma(x'')=\d^2(xc)$, etc.  satisfies $\sigma_c(I_{u/\cO})=I_{u/\cO}$. Similarly let $I^0_{u/\cO}$ be the kernel
of the map $\cO[x]\ra \cO[u]$, $x\mapsto u$, and let $\Sigma^0_{u/\cO}$  be the subgroup of $GL_n(\cO)$ consisting of all matrices $c$ such that the $\cO$-automorphism $\sigma^0_c:\cO[x]\ra \cO[x]$ defined by $\sigma_c(x)=xc$, satisfies $\sigma^0_c(I^0_{u/\cO})=I^0_{u/\cO}$.

 Here is the ``$\d$-theoretic description/bound" of the $\d$-Galois group in our cases of interest:

\begin{proposition}
\label{forget}\ 

1) $G_{u/\cO}=\Sigma_{u/\cO}$.

2)  If $\Delta(x)$ is of  type $GL_n$, $SL_n$ or $SO(q)$, we have $G_{u/\cO}\subset  G_u$. 

3) If $\Delta(x)$ is of  type $GL_n$ we have $G_{u/\cO}=\Sigma^0_{u/\cO}\cap G_u$.

\end{proposition}

 {\it Proof}. To check assertion 1 let, first, $c\in G_{u/\cO}$, so there exists $\sigma\in \tilde{G}_{u/\cO}$ with $\sigma u=c u$. Then we claim that $c\in \Sigma_{u/\cO}$. Indeed this follows from the commutativity of the diagram
\begin{equation}
\label{redd}\begin{array}{ccc}
\cO\{x\} & \stackrel{\sigma_c}{\ra} & \cO\{x\}\\
\downarrow & \  & \downarrow\\
\cO\{u\} & \stackrel{\sigma}{\ra} & \cO\{u\}
\end{array}\end{equation}
 Conversely, if $c\in \Sigma_{u/\cO}$ then $\sigma_c:\cO\{x\}\ra \cO\{x\}$ obviously induces an automorphism $\sigma:\cO\{u\}\ra \cO\{u\}$ commuting with $\d$ and sending $u$ into $uc$ so $c\in G_{u/\cO}$.
 
 To check assertions 2 and 3 we need a preliminary discussion in which we assume that $\Delta(x)$ is of type $GL_n$, $SL_n$, or $SO(q)$.
 
 Let us start with an $\cO$-algebra automorphism $\sigma$ of $\cO\{u\}$
 such that  $c:=u^{-1}\cdot \sigma(u)\in GL_n(\cO)$ and let $\epsilon=1+p\alpha$.
 We may (uniquely) extend $\sigma$ to an automorphism
 of $S^{-1}\cO\{u\}$ where $S$ is the multiplicative system consisting of all elements of $\cO\{u\}$ of the form $\det(u)^m+pf$ where $m\in \bZ_{\geq 0}$, $f\in \cO\{u\}$. 
  We claim that
\begin{equation}
\label{freeze}
\sigma(\Phi(u))=\Phi(\sigma(u)).\end{equation}
(The left hand side makes sense because $\Phi(u)=\epsilon^{-1}\cdot \phi(u)$ has entries in $S^{-1}\cO\{u\}$.)
 This is clear if $\Delta$ is of type $GL_n$ because in this case $\Phi(x)$ has polynomial entries. Let us check \ref{freeze} in case $\Delta$ is type $SO(q)$; the case of $SL_n$ is similar. Indeed, since 
 \begin{equation}
 \label{corridor}
 \Lambda(u)=(u^{(p)})^{-1}\Phi(u)=(u^{(p)})^{-1}\epsilon^{-1}\phi(u)
 =(u^{(p)})^{-1}(1+p\alpha)^{-1}(u^{(p)}+p\d u)
 \end{equation}
  it follows that $\Lambda(u)$ has entries in $S^{-1}\cO\{u\}$ so
 $$\sigma(\Phi(u))=\sigma(u^{(p)}\cdot \Lambda(u))=\sigma(u^{(p)})\cdot \sigma(\Lambda(u))=(uc)^{(p)}\cdot \sigma(\Lambda(u)),$$
 $$\Phi(\sigma u)=(uc)^{(p)}\cdot \Lambda(uc).$$
 So it is enough to check that $\sigma(\Lambda(u))=\Lambda(uc)$. Since both matrices in the latter equality are $\equiv 1$ mod $p$ in $GL_n(R)$ (for the first one use \ref{corridor}) it is enough to check that their squares are equal. But, since $M(x):=\Lambda(x)^2$ has entries rational functions of $x$, we get:
 $$\sigma(\Lambda(u))^2=\sigma(\Lambda(u)^2)=\sigma(M(u))=M(\sigma u)=
 M(uc)=\Lambda(uc)^2,$$
 which concludes the proof of \ref{freeze}.
Using \ref{freeze} in equation \ref{calculuss} below we get
\begin{equation}
\label{calculus}
\phi(\sigma(u))=\phi(uc)=\phi(u)\cdot \phi(c)=\epsilon \cdot \Phi(u) \cdot \phi(c),
\end{equation}
\begin{equation}
\label{calculuss}
\sigma(\phi(u))=\sigma(\epsilon \cdot \Phi(u))=\epsilon \cdot \sigma(\Phi(u))
=\epsilon \cdot \Phi(\sigma(u))=\epsilon\cdot \Phi(uc).
\end{equation}

To check assertion 2 let $c\in G_{u/\cO}$ and let us prove that $c\in  G_u$.
Let $\sigma\in \tilde{G}_{u/\cO}$, $\sigma u=uc$.
Since $\sigma\circ \d=\d\circ \sigma$ on $\cO\{u\}$ it follows that $\sigma\circ \phi=\phi\circ \sigma$ on $\cO\{u\}$
so, by \ref{calculus} and \ref{calculuss},  $\Phi(uc)=\Phi(u)\cdot \phi(c)$ hence $c\in G_u$.

 To check assertion 3, assume $\Delta=0$ (hence $\cO\{u\}=\cO[u]$).
 Let, first, $c\in G_{u/\cO}$ and let us prove that $c\in \Sigma^0_{u/\cO}\cap G_u$.
By assertion 2 we already know that $c\in G_u$. Also $c\in \Sigma^0_{u/\cO}$ by 
the commutativity of the diagram
\begin{equation}
\label{reddd}\begin{array}{ccc}
\cO[x] & \stackrel{\sigma^0_c}{\ra} & \cO[x]\\
\downarrow & \  & \downarrow\\
\cO[u] & \stackrel{\sigma}{\ra} & \cO[u]
\end{array}\end{equation}
Conversely let
$c\in \Sigma^0_{u/\cO}\cap G_u$  and let us prove that $c\in G_{u/\cO}$. Indeed 
since $\sigma_c^0(I^0_{u/\cO})=I^0_{u/\cO}$, it follows that $\sigma^0_c:\cO[x]\ra \cO[x]$ induces an automorphism $\sigma:\cO[u]\ra \cO[u]$  with $\sigma(u)=uc$. On the other hand since $c\in G_u$ we have
$\Phi(uc)=\Phi(u)\cdot \phi(c)$ hence, by  \ref{calculus} and \ref{calculuss}, $\phi(\sigma(u))=\sigma(\phi(u))$. It follows that $\sigma\circ \phi=\phi\circ \sigma$ on $\cO[u]$ hence $\sigma\circ \d=\d\circ \sigma$ on $\cO[u]$. So $\sigma\in \tilde{G}_{u/\cO}$, hence $c\in G_{u/\cO}$ and we are done.
\qed

\bigskip

For our discussion below we recall from the Introduction that we denote by $T, W, N$
the torus of diagonal matrices in $G$, the Weyl group of permutation matrices in $G$ and the normalizer of $T$ in $G$ respectively; so $N=TW=WT$. Also if $G^{\d}=\{a\in G;\d a=0\}$ we set $T^{\d}=T\cap G^{\d}$, $N^{\d}=N\cap G^{\d}=T^{\d}W=WT^{\d}$; $G^{\d}$ is a subset of $G$ while $T^{\d}$ and $N^{\d}$ are subgroups of $G$.

\begin{definition}
We say that $\Phi$ is right compatible with $N$ if $\Phi(ac)=\Phi(a)\cdot c^{(p)}$ for all $a\in G$ and all $c\in N$.
\end{definition}

\begin{example}
If $\Delta$ is if type  $GL_n, SL_n, SO(q)$ then $\Phi$ 
 is right compatible with $N$. By the way if $\Delta$ is of type $GL_n$  (i.e. in case $\Delta=0$) right compatibility of $\Phi(x)=x^{(p)}$ with $N$ simply means that  $(ac)^{(p)}=a^{(p)}c^{(p)}$ for $a\in G$ and $c\in N$.
\end{example}

\begin{lemma}\label{rightt}
If $\Phi$ is right compatible with $N$ then 
  $N^{\d}\subset G_u$.
\end{lemma}

{\it Proof}. Trivial.\qed
 
  \begin{lemma}\label{fuego} Assume $\Delta=0$ and set $N_{u/\cO}=N^{\d}\cap \Sigma^0_{u/\cO}$.

  1) Assume  the entries of one of the rows of $u$ are algebraically independent over $\cO$. Then $G_{u/\cO}\subset N^{\d}$ hence 
  $$G_{u/\cO}=N_{u/\cO}.$$ 
  
  2) Assume   the entries of $u$ are algebraically independent over $\cO$; then 
  $$G_{u/\cO}=N^{\d}\cap GL_n(\cO).$$
  
  3) Assume  $\sigma$ is an $\cO$-automorphism of $\cO[u]$ such that $\sigma(u)=uc$ with $c\in GL_n(\cO)\cap G_u$. Then $c\in G_{u/\cO}$.
  
  4) Assume  $n=1$. Then
$G_{u/\cO}\subset N^{\d}=G^{\d}$.

 5) We have an equality
 $$\bigcap_{u\in G}G_u=N^{\d}.$$

 \end{lemma}

{\it Proof}. 
To prove 1 let $c\in G_{u/\cO}$, hence $c\in G_u$, i.e.  $(uc)^{(p)}=u^{(p)}\phi(c)$. If $c=(c_{ij})$ then for all $m$ and $j$
 $$\sum_{i=1}^n u_{mi}^p \phi(c_{ij})=(\sum_{i=1}^n u_{mi}c_{ij})^p.$$
 Let $m$ be such that $u_{m1},...,u_{mn}$ are algebraically independent over $\cO$.
 Identifying the coefficients of the monomials in $u_{m1},...,u_{mn}$ in the latter equality we get that
  for each $j$ there exists an index $\tau(j)$ 
 such that $c_{ij}=0$ for all $i\neq \tau(j)$ and such that $c_{\tau(j)j}^p=\phi(c_{\tau(j)j})$. Since $c$ is non-singular we must have that $\tau$ is a permutation and $c\in N^{\d}$. 
 
 To prove assertion 2 note that $G_{u/\cO}\subset N^{\d}\cap GL_n(\cO)$ by assertion 1. Also   $N^{\d}\subset G_u$ by Lemma \ref{rightt}
 and, since $\cO[x]\ra \cO[u]$ is a isomorphism, we also have $\Sigma^0_{u/\cO}=GL_n(\cO)$; hence, using Proposition \ref{forget},
$N^{\d}\cap GL_n(\cO)\subset G_u\cap \Sigma^0_{u/\cO}=G_{u/\cO}$.

To prove assertion 3 let $\sigma^0_c:\cO[x]\ra \cO[x]$ be the unique $\cO$-algebra homomorphism such that  $\sigma^0_c(x)=xc$. Then $\sigma^0_c(I^0_{u/\cO})=I^0_{u/\cO}$ by the commutativity of the diagram \ref{reddd}; hence $c\in \Sigma^0_{u/\cO}$,  hence $c\in G_{u/\cO}$.
 
  To prove assertion 4 let $c\in G_{u/\cO}$; then $uc\in G^{\alpha}$ hence 
$\phi(u)\phi(c)=\epsilon u^p c^p$ where $\epsilon=1+p\alpha$. Since $\phi(u)=\epsilon u^p$ we get
 $\phi(c)=c^p$ hence $c \in G^{\d}=N^{\d}$.

To prove 5 note that the inclusion $\supset$ follows from Lemma \ref{rightt}. To prove the inclusion $\subset$
let $c$ be in the  intersection. Since $R$ is uncountable one can find $u$ with entries algebraically independent over the ring generated by the entries of $c, \d  c, \d^2 c,...$. Then one concludes that $c\in N^{\d}$ by using the same argument as in the proof of assertion 1.\qed

\begin{proposition}
\label{clock}
\

1) Assume $\Delta$ is of type $SL_n$ and  let $u\in GL_n$. Then ${\mathcal I}(x)=\det(x)$ is a prime integral for the $\Delta_u$-linear equation $\d v=\Delta_u(v)$; in other words for any $v\in G_u$
we have $\d(\det(v))=0$.

2) Assume $\Delta$ is of type $SO(q)$ and  let $u\in SO(q)$. Then ${\mathcal I}(x)=x^tqx$   is a prime integral for the $\Delta_u$-linear equation $\d v=\Delta_u(v)$; in other words for any $v\in G_u$ we have $\d(v^t q v)=0$.
\end{proposition}

{\it Proof}. To check 1) note that 
since $v\in G_u$ we have 
$$\lambda(uv)\cdot (uv)^{(p)}=\lambda(u) \cdot u^{(p)}\cdot \phi(v).$$
Taking determinants we get
$$\lambda(uv)^n\cdot \det((uv)^{(p)})=\lambda(u)^n\cdot \det(u^{(p)})\cdot \det(\phi(v)).$$
Taking into account the definition of $\lambda(x)$ we get
$$(\det(uv))^p=\det(u)^p\cdot \det(\phi(v))$$
hence $\det(v)^p=\det(\phi(v))=\phi(\det(v))$ which implies $\d(\det(v))=0$.

To check 2) note that
by Equation \ref{cry} we have
$$\phi(v)=\Phi(u)^{-1}\Phi(uv).$$
On the other hand recall that we have an identity $\Phi(x)^tq\Phi(x)=(x^t qx)^{(p)}$.
We get that 
$$\Phi(u)^tq\Phi(u)=(u^tqu)^{(p)}=q^{(p)}=q,$$
hence
$$(\Phi(u)^t)^{-1}q\Phi(u)^{-1}=q,$$
hence
$$\begin{array}{rcl}
\phi(v^tqv) & = & \phi(v)^t q \phi(v)\\
\  & \  & \  \\
\  & = & \Phi(uv)^t (\Phi(u)^t)^{-1} q \Phi(u)^{-1} \Phi(uv)\\
\  & \  & \  \\
\  & = & \Phi(uv)^t q \Phi(uv)\\
\  & \  & \  \\
\  & = & (v^t u^t q uv)^{(p)}\\
\  & \  & \  \\
\  & = & (v^tqv)^{(p)},
\end{array}
$$
which implies that $\d(v^tqv)=0$.
\qed

\medskip

Our next task will be to compute/bound the $\d$-Galois group $G_{u/\cO}$ in case $\Delta=0$. One of the morals will be that this group tends to be contained in $N$; but this is not always the case as shown in the following:

\begin{example}
\label{counterexample}
Let $\cO=\bZ_{(p)}$, $n=2$,  and assume $p\equiv 1$ mod $3$. 
Consider the matrices
$$u= \left(\begin{array}{cc} 1 & \zeta \\ 1& \zeta^2\end{array}\right),\ \ \ 
c= \left(\begin{array}{cc} 1 & $-1$ \\ 0 & $-1$\end{array}\right),\ \ \ uc= \left(\begin{array}{cc} 1 & \zeta^2 \\ 1 & \zeta\end{array}\right),$$
where $\zeta\in \bZ_p\subset R$ is a cubic root of unity.
Note that $\det u=\zeta^2-\zeta\not\equiv 0$ mod $p$ so $u, c, uc\in GL_2(R)$.
Then $u$ is a solution to the $\Delta$-linear equation
equation $$\d u=0,$$
where $\Delta=0$. We will show that $G_{u/\cO}\not\subset N$. Indeed
 $u, c, uc \in G^{\d}\backslash N$ and $u^{(p)}=u$, $c^{(p)}=c$, $(uc)^{(p)}=uc$ so $c\in G_u$.
Also we have $\cO[u]= \bZ_{(p)}[\zeta]$ and the unique non-trivial automorphism $\sigma$ of $\bZ_{(p)}[\zeta]$ sending $\sigma(\zeta)=\zeta^2$ satisfies $\sigma(u)=uc$.
By  assertion 3 in Lemma 
\ref{fuego} we have $c\in G_{u/\cO}$; so in particular
$G_{u/\cO}\not\subset N$, and our claim is proved.
  By the way
in this case 
$G_{u/\cO}=\langle c \rangle$ is cyclic of order $2$.
\end{example}

 \begin{proposition}\label{sapun}
Assume $\Delta=0$ and $\cO$ is a valuation $\d$-subring of $R$ with finite residue field. Then $G_{u/\cO}$ is a finite group.
\end{proposition}

{\it Proof}.
Let $u_0\in G^{\d}$ be the unique element such that $u\equiv u_0$ mod $p$. Let $F$ be the field of fractions of $\cO$, let $F'$ be the field generated by $F$ and the roots of unity appearing as entries in $u_0$, and let  $\cO'=R\cap F'$. Then $\cO'$ is a valuation $\d$-subring of $R$ generically finite over $\cO$ and $u_0\in GL_n(\cO')$. In particular $\cO'$ has a finite residue field.
Since $\alpha\in GL_n(\cO')$, by Proposition \ref{tata}, we get $u\in GL_n(\cO'')$ for some generically finite extension $\cO''$ of $\cO'$. Then, by the equality $\cO\{u\}=\cO[u]$, $G_{u/\cO}$ is finite.
\qed

\bigskip

 In what follows we view $R$ as a complete metric space with respect to the $p$-adic metric. So we can talk about open balls in $R$. Any open ball has the form 
$X=b+p^NR$ for some $b\in R$ and $N\in \bZ_{\geq 0}$; any such $X$ is also closed and is, again, a complete metric space with respect to the induced metric.
Now recall that a subset of a metric space is called {\it of the first category}  if it is a countable union of subsets each of which has the property that its closure has an empty interior. 
By  the Baire-Hausdorff theorem \cite{Y}, p. 11, any subset of the first category in a non-empty complete metric space $X$
is different from $X$. This applies then to any open ball $X$ in $R$.

\begin{proposition}
\label{transcendence}
Assume $\Delta=0$.
There exists  a subset $\Omega$  of the first category in the metric space
$$X=\{u\in GL_n(R);u\equiv 1\ \ \text{mod}\ \ p\}$$
such that for any $u \in X\backslash \Omega$ the following holds.
 Let $\alpha=\d u \cdot (u^{(p)})^{-1}$. Then there exists a  valuation $\d$-subring $\cO$ of $R$ containing $R^{\d}$ such that
$\alpha\in {\mathfrak g}{\mathfrak l}_n(\cO)$ and such that $G_{u/\cO}=N^{\d}$.
\end{proposition}

\begin{lemma}
\label{vanishing}
Let $x, x',..., x^{(r)}$ are a $m$-tuples of indeterminates and let $f\in R[x,x',...,x^{(r)}]\h$. Assume the 
map  $f_*:R^m\ra R$ defined by
$$f_*(a)=f(a, \d a, ..., \d^m a)$$
 vanishes on a product of open balls.  Then f vanishes on the whole of $R^m$.
\end{lemma}

{\it Proof}. By \cite{char}, Remark 1.6,  $f=0$ if and only if $f_*=0$.
So it is enough to show that for any $b_j\in R$, $1\leq j\leq m$, the  $R$-algebra homomorphism 
$$R[x,x',...,x^{(r)}]\h \ra R[x,x',...,x^{(r)}]\h,\ \ \ x_j^{(i)}\mapsto \d^i(b_j+p^Nx_j),$$
  is injective. To check this we may assume $b_j=0$ for all $j$. But then the assertion  follows from the fact that $$R[x,x',...,x^{(r)}]\h \subset K[[x,x',...,x^{(r)}]]=K[[x,\phi(x),...,\phi^r(x)]]$$ and from the fact that the endomorphism of $K[[x,\phi(x),...,\phi^r(x)]]$  defined by
$\phi^i(x)\mapsto p^N\phi^i(x)$ is injective.
\qed

\begin{lemma}
\label{Independent}
Let $E$ be a countable subfield of $K$ and let $X_1,...,X_m\subset R$ be open balls.
Then one can find 
a subset $\Omega$ of the first category in the metric space $X=X_1\times...\times X_m$ such that for all $u=(u_1,...,u_m)\in X\backslash \Omega$
 the family $$(\d^i u_j)_{i\geq 0, 1\leq j\leq m}$$ is algebraically independent over $E$.
\end{lemma}

{\it Proof}. Let  ${\mathcal F}=E[x,x',x'',...]$ be the polynomial ring
where each of $x, x', x'',...$ is an $m$-tuple of indeterminates . Hence ${\mathcal F}$ is countable. Then for each $f\in {\mathcal F}$ with $f\neq 0$  set
$$X_f:=\{u\in X;f(u,\d u,\d^2 u,...)=0\}.$$
Now we claim that each $X_f$ is closed in the metric space
$X$ and has empty interior; indeed $X_f$ is the zero locus in $X$ of  $f_*:R^m\ra R$ 
 and our claim follows from  Lemma \ref{vanishing}. The present Lemma follows now by taking $$\Omega=  \bigcup_{0\neq f\in {\mathcal F}}X_f.$$
\qed

\medskip

{\it Proof of Proposition \ref{transcendence}}.
Let $E$ be the subfield of $K$ generated over ${\mathbb Q}$ by all the roots of unity in $K$; i.e. $E={\mathbb Q}(R^{\d})$. Now $X$ in the Proposition is a product of balls so
by Lemma \ref{Independent} there exists
a subset of the first category $\Omega\subset X$ such that for all $u\in X\backslash \Omega$ 
  the family
$(\d^r u_{ij})_{r\geq0, 1\leq i,j\leq n}$
is algebraically independent over $E$. 
Let $\epsilon=\phi(u)\cdot (u^{(p)})^{-1}$, $\alpha=(\epsilon-1)/p$ and
consider the fields
$$F_s=E(\d^r \alpha_{ij}; 0\leq r\leq s,\ 1\leq i,j\leq n)=E(\phi^r (\epsilon_{ij}); 0\leq r\leq s,\ 1\leq i,j\leq n)$$
and $F=\cup_{s} F_s$. 
Let $\cO$ be a valuation $\d$-subring 
of $R\cap F$ containing $R^{\d}$ and the entries of $\alpha$ (e.g. one can take the ``maximal" choice" $\cO=R\cap F$).
Note  that for  $s\geq 1$ we have  equalities of fields
\begin{equation}
\label{fields}
E(\d^r u_{ij}; 0\leq r\leq s,\ 1\leq i,j\leq n)=F_{s-1}(u_{ij};  1\leq i,j\leq n).\end{equation}
Now the field in the left hand side of the \ref{fields} has transcendence degree $(s+1)n^2$
over $E$. Since $F_{s-1}$ has transcendence degree at most $sn^2$ over $E$
it follows from \ref{fields} that $(u_{ij})_{ij}$ are algebraically independent over $F_{s-1}$.
Since this is  true for all $s$ it follows that $(u_{ij})_{ij}$ are algebraically independent over $F$.
  By assertion 2 in Lemma \ref{fuego}, $G_{u/\cO}=N^{\d}$.
\qed

\medskip

The next Proposition shows that the $\d$-Galois group cannot be ``too large" at least if we take our data in a Zariski open set of the set of all data. In the statement below by a Zariski $K$-closed set in $GL_n(R)$ we understand the intersection of $GL_n(R)$ with a Zariski $K$-closed set of  $GL_n(K^a)$; in other words a $K$-closed set of $GL_n(R)$ is the zero set in $GL_n(R)$ of a collection  of polynomials with coefficients in $K$ in $n^2$ variables. A subgroup $\Gamma$ of $GL_n(R)$ is called diagonalizable if there exists $g\in GL_n(K^a)$ such that $g^{-1}\Gamma g$ consists of diagonal matrices.

\begin{proposition}
\label{dimless} 
There exists a Zariski $K$-closed set $\Omega$ in $G=GL_n(R)$ not containing $1$ such that for any $u\in G\backslash \Omega$ the following holds. Let $\alpha=\d u \cdot (u^{(p)})^{-1}$ and let $\cO$ be a valuation $\d$-subring of $R$ containing the entries of $\alpha$. Then $G_{u/\cO}$ contains a normal subgroup of finite index which is diagonalizable.  
\end{proposition} 

\medskip

In order to prove Proposition \ref{dimless} we need a series of Lemmas: \ref{Z}, 
\ref{later}, \ref{gee}.
In the discussion below (pertaining to these Lemmas only!)  it is convenient to temporarily change some of the notation used so far. Indeed we let ${\mathcal C}$ be an uncountable algebraically closed field of characteristic zero  (such as $K^a$ or ${\mathbb C}$) and all schemes will be
 schemes over ${\mathcal C}$. By a variety we will understand a reduced (not necessarily irreducible) scheme of finite type over ${\mathcal C}$. We use the same letter $X$ to denote 
 a variety $X$ over ${\mathcal C}$ and its set $X({\mathcal C})$ of ${\mathcal C}$-points. 
   In particular we denote by $G$ the group scheme $GL_n$ over ${\mathcal C}$ and also the ``abstract" group $GL_n({\mathcal C})$; we denote  by $T$ the  group scheme of diagonal matrices over ${\mathcal C}$ and also the ``abstract" group $T({\mathcal C})$  of diagonal matrices with entries in ${\mathcal C}$. 
   If $X$ is a variety and $x\in X$ is a point we always understand $x$ is a ${\mathcal C}$-point and we denote by $\dim_x X$ the maximum of the dimensions of the irreducible components of $X$ passing through $x$. 
   Also, in what follows,  we let $p$ be any integer $\geq 2$ (not necessarily prime).

\begin{lemma}
\label{Z}
Let $X\subset G$ be the Zariski closed subset consisting of all $v\in G$ satisfying the following properties:

1) $(v^m)^{(p)}=(v^{(p)})^m$ for all $m\geq 0$,

2) $(v^m)^{(p)}(v^{-m})^{(p)}=1$ for all $m\geq 0$.

Then $X$ has exactly one irreducible component passing through $1$ and that component is $T$.
\end{lemma}

\begin{remark}
The equalities 1) and 2) are viewed as equalities in ${\mathfrak g}={\mathfrak g}{\mathfrak l}_n({\mathcal C})$; note however that, by 1) and 2), for any $v\in X$ we have that $(v^m)^{(p)} \in G$ for all $m\in \bZ$  and hence 1) holds for all $m\in \bZ$ as an equality in $G$.\end{remark}

\begin{remark}
The set $X$ contains the group $N=WT=TW$ generated by the Weyl group $W$ and  the group $T$ of diagonal matrices with entries in ${\mathcal C}$. It is not clear whether $X$ actually coincides with the group $N$.
\end{remark}

\begin{remark}\label{sauce}
Let ${\mathbb X}$ be the closed subscheme of $G$ defined by the equations 1) and 2) in the statement of Lemma \ref{Z}; hence the variety ${\mathbb X}_{red}$
coincides  with $X$. It is interesting to note that tangent space of ${\mathbb X}$ at $1$ is the whole of the tangent space  of $G$ i.e. the Lie algebra $L(G)$ of $G$; indeed, equations 1) and 2) are easily seen to hold when $v$ is replaced by $1+\epsilon \xi$, where $\epsilon^2=0$ and $\xi$ is an arbitrary element of ${\mathfrak g}{\mathfrak l}_n({\mathcal C})$.
In particular ${\mathbb X}$ is not reduced.
\end{remark}

\medskip

{\it Proof of Lemma \ref{Z}}. 
Let $v\in X$, let $\langle v \rangle \subset G$ be the group generated by $v$, let
$H_v \subset G$ be the Zariski closure of $\langle v \rangle$ in $G$ (which is an algebraic subgroup of $G$, cf. \cite{hum}, p. 54), and let 
$H^{\circ}_v$ be the identity component of $H_v$. Clearly $H_v$ is commutative.

{\it Claim.} For all $v\in X$ we have $H_v^{\circ}\subset T$.

To check the claim note first that $\langle v \rangle \subset X$ hence $H_v\subset X$. Denote by $\Phi:G\ra {\mathfrak g}$ the map $\Phi(v)=v^{(p)}$. Clearly we have $\Phi(v^r v^s)=\Phi(v^r)\Phi(v^s)$ for all $r,s\in \bZ$ hence we have $\Phi(gh)=\Phi(g)\Phi(h)$ for all $g,h \in H_v$. Let $\varphi:H_v\ra {\mathfrak g}$ be the restriction of $\Phi$; then the regular map $\varphi$ takes values in $G$ and is a group homomorphism hence its image $H'_v:=\varphi(H_v)\subset G$ is a subgroup which is constructible. Hence $H'_v$ is a closed subgroup of $G$ (cf. \cite{hum}, p. 54) and hence $\varphi:H_v\ra H'_v$ is an algebraic group homomorphism. Consider the commutative diagram of (possibly reducible) varieties
$$\begin{array}{rcl}
H_v & \subset & G\\
\varphi \downarrow & \  & \downarrow \Phi\\
H'_v & \subset & {\mathfrak g}\end{array}$$
and the induced tangent maps between the corresponding tangent spaces at the identity
$$\begin{array}{rcl}
L(H_v) & \subset & L(G)\\
d_1\varphi \downarrow & \  & \downarrow d_1 \Phi\\
L(H'_v) & \subset & L(G)\end{array}$$
(Here $L(\ )$ denotes the Lie algebra functor. The linear map $d_1\Phi$ is  not a Lie algebra map. The map $d_1\varphi$, on the other hand, is, of course, a Lie algebra map because its source and target are abelian.) One can compute
$d_1\Phi$ explicitly: letting $v=1+\epsilon \xi\in GL_n({\mathcal C}[\epsilon])$, $\epsilon^2=0$, we have 
$$\Phi(v)=(1+\epsilon \xi)^{(p)}=diag(1+\epsilon p\xi_{11},...,1+\epsilon p\xi_{nn}).$$
 Hence the image of $d_1\Phi$ is contained in the Lie algebra $L(T)$ of the torus $T$. Since $d_1\varphi$ is surjective (because we are in characteristic zero) it follows that $L(H'_v)\subset L(T)$. Hence the identity component $(H'_v)^{\circ}$ of $H'_v$ is contained in $T$. Now, clearly,  $\Phi^{-1}(T)=T$. Hence $H_v^{\circ}\subset \Phi^{-1}((H'_v)^{\circ})\subset \Phi^{-1}(T)=T$ and our claim is proved.
 
 For any subtorus $S\subset T$ let us denote by $C(S)$ the centralizer of $S$ in $G$; moreover, for any integer $e\geq 1$ denote by $S^{1/e}$ the set of all $v\in G$ such that
 $v^e\in S$. By the above Claim and by the commutativity of $H_v$ it follows that for any $v\in X$ we have that $H_v^{\circ}$ is a subtorus of $T$ and there exists $e\geq 1$ such that
 $v\in C(H_v^{\circ})\cap H_v^{1/e}$. In particular we have
 $$X=\bigcup_{S,e} (C(S)\cap S^{1/e}\cap X)$$
 where $S$ runs through the (countable!) set of subtori  of $T$ and $e$ runs through the set of positive integers. 
 Since ${\mathcal C}$ is uncountable no irreducible variety over ${\mathcal C}$ is a countable union of proper closed subvarieties; in particular, applying this to the irreducible components of $X$ it follows 
 that there exists $e\geq 1$ and finitely many subtori $S_1,...,S_q\subset T$ such that
 \begin{equation}
 \label{deco}
 X=\bigcup_{i=1}^q (C(S_i)\cap S_i^{1/e}\cap X).\end{equation}
To conclude the proof of the Lemma we assume (as we always can) that ${\mathcal C}={\mathbb C}$. Let $V$ be an irreducible component of $X$ passing through $1$. We will prove that $V=T$ and this will end the proof. Assume $V\neq T$ and seek a contradiction. Since $V\neq T$ it follows that $V\not\subset T$ hence $V\backslash T$
is Zariski open in $V$ hence dense in $V$ in the complex topology. So there exists a sequence $x_n\ra 1$ (in the complex topology) with $x_n\in X\backslash T$.  By \ref{deco} and by replacing $x_n$ with a subsequence we may assume $x_n\in C(S_i)\cap S_i^{1/e}\cap X$ for some $i$. Let $[x_n]\in C(S_i)/S_i$ be the class of $x_n$ and choose an embedding $\rho: C(S_i)/S_i\ra GL_{\nu}({\mathcal C})$ for some $\nu$.
Then $\rho([x_n])\ra 1$ hence the eigenvalues of $\rho([x_n])$ tend to $1$.
But $[x_n]^e=1$, hence $\rho([x_n])^e=1$, for all $n$. So the eigenvalues of  $\rho([x_n])$ are $e$-th roots of unity so they form a discrete set. We get that for $n$ sufficiently big the eigenvalues of $\rho([x_n])$ are equal to $1$. But a matrix of finite order with all eigenvalues equal to $1$ must be the identity. Hence $\rho([x_n])=1$
hence $[x_n]=1$ hence $x_n\in S_i\subset T$ for some $n$, a contradiction. This ends the proof of the Lemma.
\qed

\bigskip

The next lemma is completely standard; we just include it for convenience. 

\begin{lemma}
\label{later}
Let $\pi:Z\ra Y$ be a morphism of varieties over ${\mathcal C}$ and assume $\sigma:Y\ra Z$ is a section of $\pi$. Assume $Y$ is irreducible and for $y\in Y$
consider the variety $\pi^{-1}(y)$. Let $y_0\in Y$ and assume the point $\sigma(y_0)$ is a connected component of $\pi^{-1}(y_0)$. Then there exists a Zariski open set $U\subset Y$ containing $y_0$ such that for all $y\in U$ 
the point $\sigma(y)$ is a connected component of $\pi^{-1}(y)$.
\end{lemma}

{\it Proof}.
This is a standard consequence of the semicontuinty theorem for the local dimension of fibers. Indeed let $Z^1,...,Z^m$ be the irreducible components of $Z$, let $S=\sigma(Y)$ and assume $\sigma(y_0)\in Z^i$ for $1\leq i\leq r$ and $\sigma(y_0)\not\in Z^j$ for $r< j\leq m$. Let $U_0=\pi(S\backslash \bigcup_{j>r}Z^j)$. Also let $Y^i\subset Y$ be the closure of $\pi(Z^i)$ and let $\pi_i:Z^i\ra Y^i$ for $i\leq r$ be induced by $\pi$. By the semicontinuity theorem in \cite{hum}, p.33, for $i\leq r$, there exist closed sets $T^i\subset Z^i$
not containing $\sigma(y_0)$ such that 
\begin{equation}
\label{fragrance}
\dim_x \pi_i^{-1}(\pi(x))\leq \dim_{\sigma(y_0)} \pi_i^{-1}(y_0)\ \ \text{for all}\ \ x\in Z^i\backslash T^i.\end{equation}
Consider the closed set $T:=T^1\cup...\cup T^r\cup Z^{r+1}\cup...\cup Z^m$ in $Z$
and the open subset $U=\pi(S\backslash T)=Y\backslash \pi(S\cap T)$ of $Y$. Then $y_0\in U$.
 Let $y\in U$ and let $F$ be an irreducible component of $\pi^{-1}(y)$ passing through $\sigma(y)$. Then $F\not\subset Z^j$ for $j>r$ (because if one assumes the contrary then $\sigma(y)\in S\cap Z^j\subset S\cap T$ hence $y\in \pi(S\cap T)$, a contradiction). So $F\subset Z^i$ for some $i\leq r$ and hence $F\subset \pi_i^{-1}(y)$. Since $y\not\in \pi(S\cap T)$ we have $\sigma(y)\not\in T$ hence $\sigma(y) \not\in T^i$; on the other hand $\sigma(y)\in F\subset Z^i$, hence $\sigma(y)\in Z^i\backslash T^i$. So by \ref{fragrance} we get
$$\dim_{\sigma(y)}F\leq \dim_{\sigma(y)}\pi_i^{-1}(y)\leq \dim_{\sigma(y_0)}\pi_i^{-1}(y_0)\leq \dim_{\sigma(y_0)}\pi^{-1}(y_0)=0.$$
So $\dim_{\sigma(y)}F=0$ hence $F=\{\sigma(y)\}$ and we are done.
\qed

\begin{lemma}
\label{gee} 
Let $Y$ be the Zariski open set of $G=GL_n({\mathcal C})$ consisting of all $u\in G$ such that $u^{(p)}$ is invertible. Let $\Psi:Y\times G\ra {\mathfrak g}$ be the morphism defined by 
$$\Psi(u,v) =(u^{(p)})^{-1}(uv)^{(p)}.$$
For each $u\in Y$ let
 $X_u \subset G$ be the Zariski closed set consisting of all $v\in G$ such that

1) $\Psi(u,v^m)=\Psi(u,v)^m$ for all $m\geq 0$,

2) $\Psi(u,v^m) \Psi(u,v^{-m})=1$ for all $m \geq 0$.

Then there exists a Zariski open set $U\subset Y$ containing $1$ with the property that for any $u\in U$ and for any connected closed subgroup $S\subset G$ contained in $X_u$ we have that $S$ is a torus.
\end{lemma}

{\it Proof}. Let $Z\subset Y\times G$ be the closed set defined by the equations 1) and 2) together with the equation 
 $(v-1)^n=0$.
Note that this latter equation is equivalent to asking that  $v$ be unipotent. 
Let $\pi:Z\ra Y$, $\pi(u,v)=u$,
and let $pr_G:Y\times G\ra G$ be the second projection. Then $pr_G(\pi^{-1}(u))$
coincides with the set of unipotent matrices in $X_u$. Also note that
$X_1$ coincides with $X$  in Lemma \ref{Z}. Now, by Lemma \ref{Z},
there is exactly one irreducible component of $X_1$ passing through $1$ and that component is a torus so it does not contain unipotent matrices with the exception of $1$ itself. In particular $1$ is a connected component of $pr_G(\pi^{-1}(1))$. 
 Now $\pi$ has a section $\sigma:Y\ra Z$, $\sigma(u)=(u,1)$. By Lemma \ref{later} there exists a Zariski open set $U$ of $Y$ containing $1$ such that for all $u\in U$
we have that $(u,1)$ is a connected component of $\pi^{-1}(u)$.
So $1$ is a connected component of the set of unipotent matrices in $X_u$.
Now let $S\subset G$ be a closed connected subgroup contained in $X_u$.
Then $1$ is a connected component of the set of unipotent matrices in $S$.
This implies that $S$ contains no unipotent matrix except $1$ (because any unipotent 
matrix $\neq 1$ is contained in a subgroup isomorphic to the additive group).
So the unipotent radical of $S$ is trivial,  hence a torus by \cite{hum}, p. 161.
\qed

\begin{remark}
Exactly as in Remark \ref{sauce}, if ${\mathbb X}_u$ is the subscheme of $G$ defined by equations 1) and 2) in Lemma \ref{gee} then $({\mathbb X}_u)_{red}=X_u$ and the tangent space to ${\mathbb X}_u$ at $1$ is, again, the whole of the Lie algebra $L(G)={\mathfrak g}{\mathfrak l}_n({\mathcal C})$.
\end{remark}

\medskip

{\it Proof of Proposition \ref{dimless}}. Consider the situation and notation in Lemma \ref{gee} with ${\mathcal C}=K^a$. Choose a polynomial $F\in K^a[x]$ such that
$$1\in D(F):=\{v\in GL_n(K^a);F(v)\neq 0\}\subset U.$$
Replacing $F$ by the product of its conjugates over $K$ we may assume $F\in K[x]$ and hence that $F\in R[x]$. Now let $u\in D(F)\cap GL_n(R)$, $\alpha=\d u \cdot (u^{(p)})^{-1}$, and  let $\cO\subset R$ be a valuation $\d$-subring containing the entries of $\alpha$. 
Let $\overline{G_{u/\cO}}$ be the Zariski closure of $G_{u/\cO}$ in $GL_n(K^a)$.
We want to show that the connected component $\overline{G_{u/\cO}}^{\circ}$ of $\overline{G_{u/\cO}}$ is a torus in $GL_n(K^a)$. Note that $u^{(p)}$ is invertible so $u\in Y$. Let $c\in G_{u/\cO}$ hence $c^m\in G_{u/\cO}\subset G_u$ for all $m\in \bZ$. Hence $(uc^m)^{(p)}=u^{(p)}\phi(c^m)$, hence $\Psi(u,c^m)=\phi(c^m)$.  We claim  that $c\in X_u$; indeed for $m\geq 0$ we have
$$
\Psi(u,c^m)  =  \phi(c^m)=\phi(c)^m= \Psi(u,c)^m$$
and also
$$\Psi(u,c^m)\Psi(u,c^{-m}) = \phi(c^m)\phi(c^{-m})=\phi(1)=1.$$
Since $c$ was arbitrary in $G_{u/\cO}$   we conclude that $G_{u/\cO}\subset X_u$ hence $\overline{G_{u/\cO}}\subset X_u$.
By Lemma \ref{gee}, 
 $\overline{G_{u/\cO}}^{\circ}$  is a torus. Then clearly
$$G_{u/\cO}\cap \overline{G_{u/\cO}}^{\circ}$$ is a normal subgroup of finite index in
$G_{u/\cO}$ which is diagonalizable.
\qed

\end{document}